\documentclass[reqno,a4paper]{amsart}

\usepackage{amsmath,amsfonts,amssymb}
\usepackage{verbatim}
\usepackage{enumerate}
\usepackage{url}
\usepackage{tikz}
\usepackage{setspace}
\usepackage[left=2.5cm,right=2.5cm, top=3cm, bottom=2.5cm]{geometry}
\usepackage[utf8]{inputenc} 
\usepackage[T1]{fontenc}
\usepackage{hyperref} 	
\usetikzlibrary{arrows}

\def\d{{\sf d}}

\def\E{{\sf E}}
\def\v{{\sf v}}
\def\F{{\mathcal F}}
\def\bd{{\boldsymbol{\cdot}}}

\def\N{\mathbb N}

\def\Z{\mathbb Z}

\def\ord{\mathop{\rm ord}\nolimits}

\def\lcm{\mathop{\rm lcm}}

\makeatletter
\newcommand{\dotprod}{\mathbin{\mathpalette\dotprod@\relax}}
\newcommand{\dotprod@}[2]{%
  \ooalign{$\m@th#1{\LARGE\bd}$\cr\hidewidth$\m@th#1\prod$\hidewidth\cr}%
}
\def\ldotprod{\mathop{\dotprod}\limits}
\makeatother

\theoremstyle{plain}
\newtheorem{theorem}{Theorem}[section]
\newtheorem{lemma}[theorem]{Lemma}

\newtheorem{corollary}[theorem]{Corollary}
\newtheorem{proposition}[theorem]{Proposition}

\newtheorem{example}[theorem]{Example}
\newtheorem{conjecture}[theorem]{Conjecture}
\def\proof{\noindent {\it Proof: }}

\def\qed{\hfill\hbox{$\square$}}

\theoremstyle{definition}

\newtheorem*{Ack}{Acknowledgements}
\numberwithin{equation}{section}

\title{A note on Bass' conjecture}

\author[D. V. Avelar]{D. V. Avelar}
\address{Departamento de Análise\\
Universidade Federal Fluminense\\
Niterói, RJ\\
24210-201\\
Brazil\\
}
\email{daniloavelar@id.uff.br}

\author[F. E. Brochero Mart\'{\i}nez]{F. E. Brochero Mart\'{\i}nez}
\address{Departamento de Matem\'{a}tica\\
Universidade Federal de Minas Gerais\\
Belo Horizonte, MG\\
31270-901\\
Brazil\\
}
\email{fbrocher@mat.ufmg.br }

\author[S. Ribas]{S. Ribas}
\address{Departamento de Matem\'{a}tica\\
Universidade Federal de Ouro Preto\\
Ouro Preto, MG\\
35400-000\\
Brazil\\
}
\email{savio.ribas@ufop.edu.br }

\date{\today}

\onehalfspace

\usepackage{xcolor}

\begin{document}

\maketitle

\begin{abstract}
For a finite group $G$, we denote by $\d(G)$ and by $\E(G)$, respectively, the small Davenport constant and the Gao constant of $G$. Let $C_n$ be the cyclic group of order $n$ and let $G_{m,n,s} = C_n \rtimes_s C_m$ be a metacyclic group. In \cite[Conjecture 17]{Bas}, Bass conjectured that $\d(G_{m,n,s}) = m+n-2$ and $\E(G_{m,n,s}) = mn+m+n-2$ provided $\ord_n(s) = m$. In this paper, we show that the assumption $\ord_n(s) = m$ is essential and cannot be removed. Moreover, if we suppose that  Bass' conjecture holds for $G_{m,n,s}$ and the $mn$-product-one free sequences of maximal length are well behaved, then Bass conjecture also holds for $G_{2m,2n,r}$, where $r^2 \equiv s \pmod n$.
\end{abstract}

\section{Introduction}

Let $G$ be a finite multiplicative group. A {\em sequence} $S = g_1 \bd \dots \bd g_k$ over $G$ is a finite sequence whose terms belong to $G$, the order is disregarded and repetition is allowed. 
In this sense, we may rewrite $S = \ldotprod_{g \in G} g^{[\v_g(S)]}$, where $\v_g(S)$ denotes the {\em multiplicity} of $g$ in $S$. 
The {\em length} of $S$ is $|S| = k = \sum_{g \in G} \v_g(S)$. 
We say that $T$ is a {\em subsequence} of $S$ if $\v_g(T) \le \v_g(S)$ for every $g \in G$; in this case we denote $T \mid S$ and let $S \bd T^{[-1]} = \ldotprod_{g \in G} g^{[\v_g(S) - \v_g(T)]}$. 
For a subset $K \subset G$, let $S_K = \ldotprod_{g \in K} g^{[\v_g(S)]}$ be the subsequence of $S$ formed by the terms that belong to $K$ counted with multiplicity.
For a normal subgroup $H$ of $G$, let $\overline{G} = G/H$ and let $\overline{S}$ be the sequence formed by the terms of $S$ naturally seen as elements projected in $\overline{G}$. 
We say that 
$S$ is a {\em product-one sequence} over $G$ if $1 \in \pi(S) := \{g_{\tau(1)} \dots g_{\tau(k)} \in G \mid \tau \text{ is a permutation of $[1,k]$}\}$, 
$S$ is {\em $n$-product-one sequence} if $S$ is product-one and $|S| = n$, 
$S$ is {\em product-one free} if $1 \not\in \Pi(S) := \bigcup_{T \mid S \atop |T| \ge 1} \pi(T)$, and 
$S$ is {\em $n$-product-one free} if $S$ does not contain product-one subsequences $T$ with $|T| = n$. 
The {\em small Davenport constant} $\d(G)$ of $G$ is the maximal length of a product-one free sequence over $G$, and the {\em Gao constant} $\E(G)$ of $G$ is the smallest $\ell > 0$ such that every sequence $S$ over $G$ of length $|S| \ge \ell$ has an $|G|$-product-one subsequence.

Denote by $C_n$ the cyclic group of order $n$ and let 
$$G_{m,n,s} = C_n \rtimes_s C_m = \langle x,y \mid x^m = y^n = 1, yx = xy^s \rangle$$ 
be a metacyclic group of order $mn$. For the existence of $G_{m,n,s}$, it is required that $s^m \equiv 1 \pmod n$, that is, $\ord_n(s) \mid m$, where $\ord_n(s)$ is the multiplicative order of $s$ modulo $n$.

\begin{example}\label{exemplo}
\begin{enumerate}[(a)]
\item $G_{m,n,s}$ is an abelian group if and only if $s \equiv 1 \pmod n$.

\item $G_{2,n,-1}$ is isomorphic to the well-known dihedral group $D_{2n}$ of order $2n$.

\item If $n$ is odd, then $G_{4,n,-1}$ is isomorphic to the well-known dicyclic group $Q_{4n}$ of order $4n$. In fact, $Q_{4n} \simeq \langle x, x^2y \rangle$, since $x^2y$ has order $2n$, $x^2 = (x^2y)^n$, and $(x^2y)x = x(x^2y)^{-1}$.
\end{enumerate}
\end{example}

In \cite[Lemma 4]{ZhGa}, it is proven that, for every finite group $G$,
\begin{equation}\label{zhuanggaodesig}
\E(G) \ge \d(G) + |G|.
\end{equation}
Gao \cite{Ga} proved that $\E(G) = \d(G) + |G|$ when $G$ is abelian, and Zhuang and Gao \cite{ZhGa} conjectured such equality for every finite group $G$. On the other hand, the sequence $x^{[m-1]} \bd y^{[n-1]}$ is product-one free over $G_{m,n,s}$, therefore 
\begin{equation}\label{ddesig}
\d(G_{m,n,s}) \ge m+n-2,
\end{equation} 
which implies, by Inequality \eqref{zhuanggaodesig}, that 
\begin{equation}\label{Edesig}
\E(G_{m,n,s}) \ge mn+m+n-2.
\end{equation}

Besides proving the equalities hold in Inequalities \eqref{ddesig} and \eqref{Edesig} for dihedral groups, for dicyclic groups, and for $G_{p,q,s}$ where $p, q$ are both primes, Bass \cite{Bas} proposed the following conjecture.

\begin{conjecture}[Bass {\cite[Conjecture 17]{Bas}}]\label{bassconj}
For all triples $(m,n,s) \in \N \times \N \times \Z$ such that $\ord_n(s) = m$, 
$$\d(G_{m,n,s}) = m+n-2 \quad \text{ and } \quad \E(G_{m,n,s}) = mn+m+n-2.$$
\end{conjecture}

Notice that the previous conjecture implies Zhuang \& Gao conjecture for $G_{m,n,s}$. Moreover, Bass remarked (see \cite[Section 5]{Bas}) it seems more promising to deal with $G_{m,n,s}$ satisfying $\ord_n(s) = m$ than to $\ord_n(s)$ being a proper divisor of $m$. He further explained that the proof of Conjecture \ref{bassconj} for $p, q$ primes ``{\it relies heavily on multiple elements in the same coset of the normal subgroup giving a large number of different products when multiplied in different orders. If the multiplicative order of s is large, then the methods (...) applied to $G_{m,n,s}$ will give more products, and the proof seems more likely to work.}''

Although Conjecture \ref{bassconj} has been settled for some classes of metacyclic groups such as 
\begin{enumerate}[(a)]
\item $G_{m,n,1}$, the abelian case (see \cite{ EGZ} and Inequality \eqref{ddesig} for $m=1$ and \cite{Ols2} for $m > 1$),
\item $G_{m,mn,s}$ (see \cite{HZ}), 
\item $G_{p,n,s}$, where $p$ is the smallest prime divisor of $pn$ and $\gcd(p(s-1), n) = 1$ (see \cite{QL}), and 
\item $G_{2,n,s}$ (see \cite{Bas} for $s\equiv-1 \pmod n$ and \cite{AMR} for $s \not\equiv\pm1 \pmod n$), 
\end{enumerate} 
in this paper we prove that the hypothesis $\ord_n(s) = m$ cannot be removed in general. We deal with the case $m \ge 4$ even, $n \ge 3$ odd and $\ord_n(s) = 2$. Yet, we believe that Zhuang \& Gao conjecture holds. It is worth mentioning that two examples of small groups $G_{m,n,s}$ that do not fit $\ord_n(s) = m$ and do not satisfy Conjecture \ref{bassconj} have already been found computationally in \cite[Section 2]{CDS}. In this paper, we also prove that if Bass' conjecture holds for $G_{m,n,s}$ and the $mn$-product-one free sequences of maximal length are well behaved, then the same holds for $G_{2m,2n,r}$, where $r^2 \equiv s \pmod n$.  

The paper is organized as follows. In Section \ref{secproofnodd}, we assume $m$ even, $n$ odd and $\ord_n(s) = 2$ in order to prove that hypothesis $\ord_n(s) = m$ is essential in Conjecture \ref{bassconj}. We further provide infinitely many concrete examples for the case $m=4$. In Section \ref{secthmneven}, for any integers $m,n$, we prove that if Bass' conjecture holds for $G_{m,n,s}$ and the $mn$-product-one free sequences over $G_{m,n,s}$ of maximal length are well behaved in some expected sense, then Bass' conjecture also holds for $G_{2m,2n,r}$, where $r^2 \equiv s \pmod n$.

\section{The hypothesis $\ord_n(s) = m$ is essential in Bass' conjecture}\label{secproofnodd}

In this section, inspired by the isomorphism given in Example \ref{exemplo}(c), we prove that in general the assumption $\ord_n(s) = m$ cannot be removed in Conjecture \ref{bassconj}. We start with an upper bound for the small Davenport constant.

\begin{lemma}[See \cite{OW}]\label{lemcotad}
For any finite non-cyclic group $G$, $\d(G) \le \lceil |G|/2 \rceil$.
\end{lemma}

The following result provides lower bounds that are often higher than those given by Inequality \ref{ddesig} and \ref{Edesig}.

\begin{proposition}\label{thmhypessential}
Let $m, n$ be positive integers, and let $s$ be an integer such that $\ord_n(s) \mid m$ and $s \not\equiv 1 \pmod n$. Suppose that $G_{m,n,s}$ contains a cyclic normal subgroup $H = \langle z \mid z^t=1 \rangle$ of order $t$. Then 
$$\d(G_{m,n,s}) \ge t \quad \text{ and } \quad \E(G_{m,n,s}) \ge mn + t.$$
In particular, if $t \ge m+n-1$, then $G_{m,n,s}$ does not satisfy the conclusion of Conjecture \ref{bassconj}.
\end{proposition}

\proof
Since $s \not\equiv 1 \pmod n$, $G_{m,n,s}$ is a non-abelian group, thus $H$ is a proper subgroup of $G_{m,n,s}$. Let $w \in G_{m,n,s} \backslash H$.
The sequence $z^{[t-1]} \bd w$ is product-one free over $G_{m,n,s}$, which implies that $\d(G_{m,n,s}) \ge t$. The second inequality follows from Inequality \eqref{zhuanggaodesig}.

\qed

\begin{corollary}\label{corollary}
Let $m \ge 4$ be even, $n \ge 3$ be odd, and $\ord_n(s) = 2$. Assume that $\lcm(\frac{m}{2},n) \ge m+n-1$. Then $G_{m,n,s}$ does not satisfy the conclusion of Conjecture \ref{bassconj}. In particular, if $\gcd(\frac{m}{2},n) = 1$, then $\d(G_{m,n,s}) = \frac{mn}{2}$, $\E(G_{m,n,s}) \ge \frac{3mn}{2}$ and the conclusion of Conjecture \ref{bassconj} does not hold.
\end{corollary}

\proof
The group $H = \langle x^2y \rangle$ is a cyclic normal subgroup of $G_{m,n,s}$ of order $\lcm(\frac{m}{2},n)$, therefore the first part of the statement follows directly from the previous proposition. For the second part, since $\lcm(\frac{m}{2},n) = \frac{mn}{2} \ge m+n-1$ for every $m \ge 4$ and $n \ge 3$, the result follows from the previous proposition using that $\d(G_{m,n,s}) \le \frac{mn}{2}$ by Lemma \ref{lemcotad}.

\qed

\begin{example}
Let $m=4$ and let $n = p^t$ for some prime $p \equiv 3 \pmod 4$, where $t \ge 1$ is an integer. 
Note that $\ord_n(s) \neq 4$, for otherwise $s^2 \equiv -1 \pmod p$ which is impossible when $p \equiv 3 \pmod 4$.
Moreover, if $s^4 \equiv 1 \pmod n$, then either $s \equiv 1 \pmod n$ or $s \equiv -1 \pmod n$. Example \ref{exemplo}(c) ensures that the non-abelian group $G_{4,p^t,-1}$ is isomorphic to $Q_{4n}$. 
As there are infinitely many primes $p$ with $p \equiv 3 \pmod 4$, by Corollary \ref{corollary} and together with the above discussion, there are infinitely many groups of the form $G_{4,p^t,-1}$ (note that $\gcd(2,p^t) = 1$ and $\lcm(2,p^t) \ge 4 + p^t - 1$) for which $\d(G_{4,p^t,-1}) > 4 + p^t - 2$. Therefore, the hypothesis $\ord_n(s) = m$  is essential in Bass' conjecture.
The examples computationally found in \cite[Section 2]{CDS} are of this form, namely $G_{4,3,-1}$ and $G_{4,7,-1}$.
\end{example}

\section{Bass' conjecture for $G_{2m,2n,s}$ through $G_{m,n,s^2}$}\label{secthmneven}

In this section, we consider the even integers $m \ge 4$ and $n \ge 6$, say $m = 2m_0$ and $n = 2n_0$, and we let $s$ be an odd integer such that $\ord_{n_0}(s^2) = m_0$. There is no obstruction in assuming $s$ is odd. In fact, if $n_0$ is even, then $s$ is automatically odd, and if $n_0$ is odd and $s$ is even, then we consider the odd number $s+n_0$ modulo $n$. Moreover, since $s$ is odd, it follows that $\ord_{n_0}(s^2) = \ord_n(s)/2$, in such way that $\ord_{n_0}(s^2) = m_0$ if and only if $\ord_n(s) = m$. However, it is not really needed.

It is worth mentioning that $G_{m,n,s}$ do always have $G_{m_0,n_0,s^2}$ as a normal subgroup, and $G_{m_0,n_0,s^2}$ is a normal subgroup of both $G_{m,n,s}$ and $G_{m,n,n-s}$, which are non-isomorphic groups (actually, instead of $n-s$ one could have any other square root of $s^2$ modulo $n$). On the other hand, if $r_0 \in [0,n-1]$ is odd, $\ord_{n_0}(r_0) = m_0$ and $r_0$ is not a quadratic residue modulo $n$, then $G_{m_0,n_0,r_0}$ is not a normal subgroup of any group of the form $G_{m,n,r}$.

Consider the following somewhat natural conditions.
\begin{enumerate}[(A)]
\item $\E(G_{m_0,n_0,s^2}) = m_0n_0 + m_0 + n_0 - 2$.
\item If $S \in \F(G_{m_0,n_0,s^2})$ has length $|S| = \E(G_{m_0,n_0,s^2}) - 1$ and is $m_0n_0$-product-one free, then 
$$S = (y^{\alpha})^{[\ell n_0-1]} \bd (y^{\beta})^{[(m_0 + 1 - \ell)n_0-1]} \bd \ldotprod_{1 \le i \le m_0-1} (x^wy^{\gamma_i}),$$ 
where $\alpha, \beta, \gamma_i \in [0,n_0-1]$ with $\gcd(\alpha-\beta,n_0)=1$, $w \in [1,m_0-1]$ with $\gcd(w,m_0)=1$, and $\ell \in [1,m_0]$.
\end{enumerate}

When $\ord_{n_0}(s^2) = m_0$, the condition (A) is equivalent to Conjecture \ref{bassconj}. We further observe that the condition (B) holds for almost all metacyclic groups $G$ for which the $|G|$-product-one free sequences $S \in \F(G)$ of maximal length $|S| = \E(G)-1$ have been found so far (see \cite{AMR, OhZh, QL2}), but there is at least one genuine exception (see \cite[Theorem 1.2]{OhZh}). Indeed the group $G_{2,3,-1}$ contains an extra $6$-product-one free sequence of length $8$, namely $S = 1^{[5]} \bd x \bd xy \bd xy^2$.

In this section, we prove that conditions (A) and (B) imply $\E(G_{m,n,s}) = mn + m + n - 2$. We use the quotient group $G_{m,n,s}/G_{m_0,n_0,s^2} = \overline{G_{m,n,s}} \simeq C_2 \times C_2$ and explore the subproducts in $G_{m_0,n_0,s^2}$ and the remainder terms in $C_2 \times C_2$. The main result of this paper is the following.

\begin{theorem}\label{thmneven}
Let $m = 2m_0 \ge 4$, $n = 2n_0 \ge 6$ be integers, and let $s$ be an odd integer such that $\ord_{n_0}(s^2) = m_0$. If conditions (A) and (B) hold for $G_{m_0,n_0,s^2}$, then $\E(G_{m,n,s}) = mn + m + n -2$. 
\end{theorem}

\proof
By Inequalities \eqref{zhuanggaodesig}, \eqref{ddesig} and \eqref{Edesig}, it is only required to obtain a tight upper bound for $\E(G_{m,n,s})$. In this sense, let $S$ be a sequence over $G_{m,n,s}$ of length $|S| = mn+m+n-2 = 4m_0n_0 + 2m_0 + 2n_0 - 2$.

Let $H = \langle x^2, y^2 \rangle \simeq C_{n_0} \rtimes_{s^2} C_{m_0}$ be a normal subgroup of $G_{m,n,s}$. We have $\overline{G_{m,n,s}} \simeq C_2 \times C_2$. By Pigeonhole Principle, each $5$ terms in $G_{m,n,s}$ yield $2$ in the same class into $\overline{G_{m,n,s}}$, whose product belongs to $H$. Thus it is possible to obtain disjoint subsequences $T_1, \dots, T_{2m_0n_0 + m_0 + n_0 - 3} \mid S$, each of length $|T_i| = 2$, such that $\pi(T_i) \subset H$. Let $h_i \in \pi(T_i)$ for $1 \le i \le 2m_0n_0 + m_0 + n_0 - 3$. Since $2m_0n_0 + m_0 + n_0 - 3 > m_0n_0 + m_0 + n_0 - 2 = \E(G_{m_0,n_0,s^2})$ and $H \simeq C_{n_0} \rtimes_{s^2} C_{m_0}$, condition (A) ensures that $h_1 \bd \dots \bd h_{2m_0n_0 + m_0 + n_0 - 3}$ contains an $m_0n_0$-product-one subsequence, say $h_{m_0n_0 + m_0 + n_0 - 2} \bd \dots \bd h_{2m_0n_0 + m_0 + n_0 - 3}$. It yields a $2m_0n_0$-product-one subsequence of $S$. We now look at the sequence $h_1 \bd \dots \bd h_{m_0n_0 + m_0 + n_0 - 3}$. If it contains an $m_0n_0$-product-one subsequence, then we are done. Otherwise, condition (B) ensures that such sequence is $m_0n_0$-product-one free over $H$ if and only if
$$h_1 \bd \dots \bd h_{m_0n_0 + m_0 + n_0 - 3} = (y^{2a})^{[\ell n_0-1]} \bd (y^{2b})^{[(m_0+1-\ell) n_0-1]} \bd \ldotprod_{1 \le j \le m_0-1} x^{2d}y^{2c_j},$$
where $\gcd(a-b,n_0) = 1$, $\gcd(d,m_0) = 1$ and $\ell \in [1,m_0]$. We may assume without loss of generality that $\ell \le \lfloor \frac{m_0+1}{2} \rfloor$. Since $(y^{2a})^{n_0} = 1 = (y^{2b})^{n_0}$, it is possible to extract several $n$-product-one subsequences of $T_1 \bd \dots \bd T_{m_0n_0 + m_0 + n_0 - 3}$. We reindex if needed and set
\begin{equation}\label{hs}
h_j = 
\begin{cases}
y^{2a} &\text{ for } j \in [1, n_0-1], \\
y^{2b} &\text{ for } j \in [n_0, 3n_0-2], \\
x^{2d}y^{2c_j} &\text{ for } j \in [3n_0-1, 3n_0+m_0-3].
\end{cases}
\end{equation}
We also assume that 
\begin{equation}\label{piTs}
1 \in \pi(T_{3n_0+m_0-2} \bd \dots \bd T_{m_0n_0 + m_0 + n_0 - 3}). 
\end{equation}
Therefore, we just need to prove that $T_1 \bd \dots \bd T_{3n_0+m_0-3} \bd R$ contains a $2n$-product-one subsequence, where $R = S \bd (T_1 \bd \dots \bd T_{2m_0n_0+m_0+n_0-3})^{[-1]}$, $|R| = 4$ and $\overline{R}$ is $2$-product-one free over $\overline{G_{m,n,s}}$. The sequence given by Equation \eqref{hs} implies that $\pi(T_i) = \{h_i\}$ for $1 \le i \le 3n_0-2$, otherwise we could change the order of one of the products and obtain a distinct sequence.

Let $R_1 = T_1 \bd \dots \bd T_{3n_0+m_0-3}$ and $R_2 = T_{3n_0 + m_0 - 2} \bd \dots \bd T_{2m_0n_0 + m_0 + n_0 - 3}$, so that $1 \in \pi(R_2)$ and $|R_2| = mn - 2n$. Notice that $R = S \bd (R_1 \bd R_2)^{[-1]}$, $R$ has $4$ terms and no two of them should result in a product in $H$, otherwise we are done. Therefore 
$$R = (x^{u_1}y^{v_1}) \bd (x^{u_2}y^{v_2}) \bd (x^{u_3}y^{v_3}) \bd (x^{u_4}y^{v_4}),$$ 
where $u_i \in [0,m-1]$ and $v_i \in [0,n-1]$ for $1 \le i \le m$, $u_1, u_2, v_1, v_3$ are even and $u_3, u_4, v_2, v_4$ are odd. Let $x^{2e}y^{2f} \in \pi(R) \subset H$ and let $R_3 = x^{2e}y^{2f} \bd \ldotprod_{1 \le j \le m_0-1} x^{2d}y^{2c_j} \in \F(H)$. Notice that $|R_3| = m_0 = \d(C_{m_0}) + 1 = \d(H/C_{n_0}) + 1$. Thus $R_3$ contains a subsequence $T$ of length $1 \le |T| \le m_0$ such that $\pi(T) \subset \langle y^2 \rangle \simeq C_{n_0}$. We observe that $|T|$ terms from $R_3$ correspond to either $r = 2|T|$ or $r = 2|T|+2$ terms of $S$, depending on whether $x^{2e}y^{2f} \mid T$ or not, respectively. Let us fix $y^{2t} \in \pi(T)$ and $k \in [1,n_0]$ such that $k \equiv (b-a)^{-1}(b \cdot \frac{r}{2} - t) \pmod {n_0}$,  which actually exists since $\gcd(a-b , n_0) = 1$ and $r$ is even. 
We observe that $n_0-k\in[0,n_0-1]$ and $n_0+k- \frac r2 \in [n_0-m_0,2n_0-1] \subset [1,2n_0-1]$ once $m \mid \varphi(n) < n$, where $\varphi$ denotes the Euler totient function. We have
$$(y^{2a})^{n_0-k} \cdot (y^{2b})^{n_0+k-\frac r2} \cdot y^{2t} = 1.$$
Therefore $R_1 \bd R = S \bd R_2^{[-1]}$ contains a $2n$-product-one subsequence. Since $R_2$ is an $(mn-2n)$-product-one sequence, it follows that $S$ contains an $mn$-product-one subsequence. This implies that $\E(G_{m,n,s}) = mn + m + n - 2$.
\qed

\Ack The authors would like to thank Steven J. Miller for the valuable discussion that improved the presentation of Theorem 3.1. F.E. Brochero Martínez and S. Ribas were partially supported by FAPEMIG grants RED-00133-21 and APQ-02546-21, Brazil. F.E. Brochero Martínez was also partially supported by FAPEMIG grant APQ-02973-17, Brazil.

\end{document}